# Axiomatic (and Non–Axiomatic) Mathematics

**Saeed Salehi**

15 August 2020

Table 1: Axiomatizability

|       | $\mathbb{N}$ | $\mathbb{Z}$ | $\mathbb{Q}$ | $\mathbb{R}$ | $\mathbb{C}$ |
|-------|---|---|---|---|---|
| $\{<\}$ | ✓ | ✓ | ✓ | ✓ | – |
| $\{+\}$ | ✓ | ✓ | ✓ | ✓ | ✓ |
| $\{<,+\}$ | ✓ | ✓ | ✓ | ✓ | – |
| $\{+,\times\}$ | × | × | × | ✓ | ✓ |
| $\{\times\}$ | ✓ | ✓ | ✓ | ✓ | ✓ |
| $\{<,\times\}$ | × | × | ✓ | ✓ | – |
| exp | × | – | – | ? | × |

Table 2: Axioms for Identities (over $\mathbb{R}^+$)

| | |
|---|---|
| $\{+\}$ | $x+(y+z)=(x+y)+z,\quad x+y=y+x$ |
| $\{\times\}$ | $x\cdot(y\cdot z)=(x\cdot y)\cdot z,\quad x\cdot y=y\cdot x,\quad x\cdot\mathbf{1}=x$ |
| $\{+,\times\}$ | $x\cdot(y+z)=(x\cdot y)+(x\cdot z)$ |
| $\{\text{exp}\}$ | $(x^y)^z=(x^z)^y,\quad x^1=x,\quad \mathbf{1}^x=\mathbf{1}$ |
| $\{\times,\text{exp}\}$ | $x^{(y\cdot z)}=(x^y)^z,\quad (x\cdot y)^z = x^z\cdot y^z$ |
| $\{+,\times,\text{exp}\}$ | $x^{(y+z)}=x^y\cdot x^z,\quad \cdots$ |

*The method of "postulating" what we want has many advantages;
they are the same as the advantages of theft over honest toil.*
–Bertrand Russell (1919, *Introduction to Mathematical Philosophy*)

---

**A**XIOMATIZING mathematical structures and theories, or *postulating* them as Russell put it, is an objective of Mathematical Logic. Some axiomatic systems are nowadays mere definitions, such as the axioms of *Group Theory*; but some systems are much deeper, such as the axioms of *Complete Ordered Fields* with which Real Analysis starts. Groups abound in mathematical sciences, while by Dedekind's theorem there exists only one complete ordered field, up to isomorphism. Cayley's theorem in Abstract Algebra implies that the axioms of group theory completely axiomatize the class of permutation sets that are closed under composition and inversion; cf. e.g. [7].

In this article, we survey some old and new results on the first-order axiomatizability of various mathematical structures (Table 1). The (non-)axiomatizability of many structures in Table 1 are known from almost a century ago; for example, the axiomatizability of $\langle\mathbb{C};+,\times\rangle$ follows from Tarski's theorem (1936), and the non-axiomatizability of $\langle\mathbb{N};+,\times\rangle$ follows from Gödel's theorem (1931). The question of the axiomatizability of e.g. $\langle\mathbb{Q};<,\times\rangle$ seemed to be missing in the literature, which was shown to be axiomaitzable in [1] for the first time; Tarski's result implies the axiomatizability of $\langle\mathbb{C};\times\rangle$, but one explicit axiomatization for it was presented in [20] for the first time. We will also review identities over $+,\times,\text{exp}$ that hold in the set of positive real numbers (Table 2). The identities on Table 2, except the last row which contains three dots, do completely axiomatize the identities that hold in the set of positive real numbers ($\mathbb{R}^+$) over the indicated operations. Whether all the identities in the table completely axiomatize the identities in the structure $\langle\mathbb{R}^+;\mathbf{1},+,\times,\text{exp}\rangle$ is the well-known *Tarski's High-School Problem*, which has an interesting history.

## § Boolean Algebras & Propositional Logic.

Arguably, Modern Logic starts with Boole's *Investigation of the Laws of Thought* (1854); Boole's axiomatic system is called "propositional logic" nowadays. It axiomatizes some of the basic properties of the conjunction (∧), disjunction (∨), and negation (¬) connectives. The Boolean expressions (or Boolean terms, or propositional formulas) are constructed from a fixed infinite set of atoms, say $\{p_0,p_1,p_2,\cdots\}$, by means of those connectives. Let us note that implication (→) is definable by disjunction and negation as $(a\to b)\equiv(\neg a)\vee b$, where $\equiv$ denotes (logical) equivalence. Boole's axiomatization is in fact nothing but a definition of Boolean Algebras:

---

Associativity
  $a\wedge(b\wedge c)\equiv(a\wedge b)\wedge c,\quad a\vee(b\vee c)\equiv(a\vee b)\vee c$
Commutativity
  $a\wedge b\equiv b\wedge a,\quad a\vee b\equiv b\vee a$
Distributivity
  $a\wedge(b\vee c)\equiv(a\wedge b)\vee(a\wedge c),\quad a\vee(b\wedge c)\equiv(a\vee b)\wedge(a\vee c)$
Idempotence
  $a\wedge a\equiv a,\quad a\vee a\equiv a$
Truth and Falsum
  $a\vee(\neg a)\equiv\top,\quad a\wedge\top\equiv a,\quad a\wedge(\neg a)\equiv\bot,\quad a\vee\bot\equiv a$
de Morgan's Laws
  $\neg(a\wedge b)\equiv(\neg a)\vee(\neg b),\quad \neg(a\vee b)\equiv(\neg a)\wedge(\neg b)$

---

Many more identities can be deduced (proved) from the above axioms, such as the following:

**EXAMPLE 1** (i) It immediately follows from the axioms that $a \equiv a \wedge \top \equiv a \wedge (p \vee \neg p) \equiv (a \wedge p) \vee (a \wedge \neg p)$.

(ii) The *absorbing properties* of truth and falsum, i.e., $a \vee \top \equiv \top$ and $a \wedge \bot \equiv \bot$ follow also from the axioms. We show the former: $a \vee \top \equiv a \vee (a \vee \neg a) \equiv (a \vee a) \vee (\neg a) \equiv a \vee (\neg a) \equiv \top$.

(iii) One can also prove *the absorption laws*: $a \wedge (a \vee b) \equiv a$ and $a \vee (a \wedge b) \equiv a$. Let us show the latter by using (ii) above: $a \vee (a \wedge b) \equiv (a \wedge \top) \vee (a \wedge b) \equiv a \wedge (\top \vee b) \equiv a \wedge (b \vee \top) \equiv a \wedge \top \equiv a$.

(iv) The *double negation law* $\neg \neg a \equiv a$ can be proved as follows: $\neg \neg a \equiv (\neg \neg a) \wedge \top \equiv (\neg \neg a) \wedge (a \vee \neg a) \equiv (\neg \neg a \wedge a) \vee (\neg \neg a \wedge \neg a) \equiv (\neg \neg a \wedge a) \vee (\bot) \equiv (a \wedge \neg \neg a) \vee (a \wedge \neg a) \equiv a \wedge (\neg \neg a \vee \neg a) \equiv a \wedge \top \equiv a$. ◇

We show that *all* the valid laws, according to the truth-table semantics, are provable from the axioms; thus it is a complete axiomatic system (for Boolean equivalences):

**THEOREM 2 (Completeness)** *If $a \equiv b$ is valid according to the truth-table semantics, then it is provable from the axioms.*

A proof can proceed by normalizing the Boolean terms, or propositional formulas. A (propositional) formula $a$ is said to be in *disjunctive normal form* (DNF) when it is a disjunction of some formulas each of which is a conjunction of some atoms or negated atoms; i.e., $a = \bigvee_i c_i$ where each $c_i$ is $\bigwedge_j \ell_{(i,j)}$ for some atoms or negated-atoms $\ell_{(i,j)}$.

If $p$ is an atom, then $(p)$ and $(\neg p)$ are both DNF; if $q$ is another atom, then the four formulas $(p) \vee (q)$, $(p) \vee (\neg p \wedge q)$, $(p \wedge \neg q) \vee (q)$, and $(p \wedge q) \vee (p \wedge \neg q) \vee (\neg p \wedge q)$ are equivalent DNF's. Every propositional formula can be seen to be equivalent to a DNF formula, and this can be proved by the above axioms: firstly implication ($\rightarrow$) does not appear in our formulas; and secondly by the double negation law, proved in Example 1(iv), and de Morgan's laws, negations ($\neg$) can be pushed as far as possible inside the sub-formulas, so that they appear at most behind atoms. Finally, by distributing all the conjunctions over disjunctions, if any, an equivalent DNF formula is obtained; and this equivalence is provable from the above axioms.

Now the proof goes as follows: assume that all the atoms that appear in $a$ and $b$ belong to the set $\{p_0, \cdots, p_k\}$; $a$ and $b$ are provably equivalent to some DNF formulas, such as e.g. $a \equiv \bigvee_i c_i$ and $b \equiv \bigvee_j d_j$ where $c_i$'s and $d_j$'s are conjunctions of some atoms or negated atoms. By Example 1(i) we can assume that all the atoms $p_0, \cdots, p_k$ appear exactly once in each $c_i$ and $d_j$. By this assumption, we show that each $c_i$ is equal to some $d_j$, and vice versa. Thus, $a$ and $b$ are provably equivalent. For a fixed $c_i$ consider the evaluation that maps an atom to $\top$ if it appears positively in $c_i$, and maps it to $\bot$ if it appears negatively in $c_i$. Under that evaluation, $c_i$, and so $a$, is mapped to $\top$; thus $b$ should be mapped to $\top$ too. So, some $d_j$ should be mapped to $\top$ under that evaluation; and this is possible only when $d_j = c_i$.

The completeness of Propositional Logic with respect to truth-table semantics follows from Theorem 2. For example, the validity of the formula $[(p \rightarrow q) \rightarrow p] \rightarrow p$, Peirce's Law, can be proved by first translating $a \rightarrow b$ to $\neg a \vee b$, and then showing the equivalence $(\neg[\neg(\neg p \vee q) \vee p] \vee p) \equiv \top$ by the above axioms.

We will come back to mathematical *Identities* at the end of the paper; before that let us study the axiomatizabilitiy of some mathematical structures.

## § AXIOMATIZABILITY & QUANTIFIER ELIMINATION.

A first-order structure consists of a non-empty set $D$, which is called domain (universe), together with a first-order language $\mathcal{L}$, consisting of some constant, relation or function symbols that are interpreted over the domain. The abstract definition of a structure $\mathfrak{A} = \langle D; \mathcal{L} \rangle$ from Model Theory is not needed here (see e.g. [17] for more details). In the first-order setting, the quantifiers ($\forall, \exists$) range over the elements of the domain in question (which are taken to be number sets $\mathbb{N}, \mathbb{Z}, \mathbb{Q}, \mathbb{R}$, and $\mathbb{C}$, here). So, subsets of the domain cannot be quantified; thus, the statement "for every nonempty and bounded subset there is a supremum for it" is not first-order, while "every element has an inverse" is so.

One reason for studying mathematical structures and theories in the setting of first-order logic is that despite of the fact that this logic is too weak to represent some fundamental properties (such as *begin well-ordered* or *completeness* of ordered sets) it has some other nice properties such as the *compactness* and *semantic completeness* (proved by Gödel 1930).

On the other hand, second-order logic may seem to be a more expressive framework for studying mathematical theories and structures (in which one can express the properties of *begin well-ordered* and *completeness* of ordered sets). But it has its own foundational problems; the same problems that set theory has with *incompleteness* and *truth* (proved by Gödel 1931). In fact, as Quine put it, the second-order logic is "set theory in sheep's clothing" (this is actually the title of the fourth section of the fifth chapter of Quine's *Philosophy of Logic*, 1986).

So, we have chosen first-order logic as the framework of our study; though, the study could be undertaken in the framework of second-order logic as well. Let us recall that a *sentence* is a formula without any free variables, i.e., all of its variables are quantified; and a *theory* is a set of sentences. We saw in the previous section that propositional logic is axiomatizable; so a way of axiomatizing a structures is reducing its first-order theory to propositional logic which is usually done through the process of Quantifier Elimination.

**DEFINITION 3 (Quantifier Elimination, QE)** A theory $T$ is said to admit *quantifier elimination* (QE) when there exists an algorithm that for a given formula $\varphi(\vec{x})$ as input, with the shown free variables, outputs a quantifier-free formula $\theta(\vec{x})$ with exactly the same free variables ($\vec{x}$) such that $T$ proves the sentence $\forall \vec{x}[\varphi(\vec{x}) \leftrightarrow \theta(\vec{x})]$. ◇

So, if $T$ admits QE, then every first-order sentence over its language is equivalent in $T$ to an algorithmically calculable quantifier-free sentence. Quantifier elimination is usually done by the means of the following fundamental lemma which is proved also in [6, Theorem 31F], [13, Theorem 4.1], and [23, Lemma III.4.1].

**LEMMA 4 (The Main Lemma of Quantifier Elimination)**
*A theory $T$ admits QE if and only if there exists an algorithm that for every given formula of the form $\exists x \gamma(x)$, where $\gamma(x)$ is a conjunction of some atoms or negated atoms, outputs a quantifier-free formula $\theta$ such that the free variables of $\theta$ are all the free variables of $\gamma(x)$ other than $x$, and the universal closure of $[\exists x \gamma(x) \leftrightarrow \theta]$ is provable in $T$.*

The "only if" part of the lemma is trivial. For the "if" part, let $\varphi$ be an arbitrary formula. We show that it is $T$–equivalent to a quantifier-free formula with the same free variables (as of $\varphi$) and that quantifier-free formula can be found algorithmically. Take one of the innermost quantifiers of $\varphi$; such as $\forall x\theta(x)$ or $\exists x\theta(x)$ where $\theta$ is a quantifier-free formula. In the former case consider $\neg\exists x\neg\theta(x)$; so without loss of generality we can assume that the quantifier is existential. We saw (in the Proof of Theorem 2) that every propositional formula is equivalent to a DNF formula. So, $\exists x\theta(x) \equiv \exists x \bigvee_i \gamma_i(x) \equiv \bigvee_i \exists x \gamma_i(x)$, where each $\gamma_i(x)$ is a conjunctions of some atomic or negated atomic formulas. By the assumption, the existing algorithm can find a $T$–equivalent quantifier-free formula for each $\exists x \gamma_i(x)$; thus that algorithm can find a $T$–equivalent formula for $\varphi$ with one less quantifier (than $\varphi$). So, by an inductive argument one can show the existence of an algorithm that outputs a quantifier-free formula with the same free variables (as of $\varphi$) that is moreover $T$–equivalent to $\varphi$.

Quantifier Elimination is applicable for axiomatizing the complete first-order theory of a structure $\mathfrak{A}$ when we have a candidate theory $T$ in a way that (i) all the axioms of $T$ are true in $\mathfrak{A}$, (ii) $T$ admits QE, and (iii) $T$ decides (i.e., either proves or refutes) every atomic sentence. Then, $T$ is a *complete* theory, in the sense that it either proves or refutes every sentence over the language of $T$, and so it *completely axiomatizes* $\mathfrak{A}$. Thus, $T$ proves every sentence that is true in $\mathfrak{A}$, and refutes every sentence that is not true in $\mathfrak{A}$. In the following, we will study some axiomatizations of number systems ($\mathbb{N}, \mathbb{Z}, \mathbb{Q}, \mathbb{R}, \mathbb{C}$) over the first-order languages that may contain $<$, $+$, $\times$, or $\exp$.

For a structure that is known to be axiomatizable (we do not have a clear criterion for axiomatizability or non-axiomatizability of a given structure), we introduce a theory that is true in that structure and decides every atomic sentence, and show that it admits QE. Thus, the proposed theory does completely axiomatize the structure.

### § Number Systems (Order & Addition).

Let us first study the order relation ($<$) in number systems. We recall that an *order* is a binary relation that is anti-symmetric, transitive, and linear (see the axioms **A**$_<$, **T**$_<$, **L**$_<$ in Theorem 5). The order is dense in $\mathbb{Q}$ and $\mathbb{R}$ (see **D**$_<$ in Theorem 5) and has no endpoints (see **U**$_<$ and **B**$_<$ in Theorem 5). This is all the first-order theory of order can say in $\mathbb{Q}$ and $\mathbb{R}$, since it is a complete theory. However, the structure $\langle \mathbb{Q}; < \rangle$ is very different from the structure $\langle \mathbb{R}; < \rangle$, since the latter is complete (every nonempty and bounded subset has a supremum) while the former is not.

**THEOREM 5 (An Axiomatization for $\langle \mathbb{R}; < \rangle$ and $\langle \mathbb{Q}; < \rangle$)**
*The (finite) theory of dense linear orders without endpoints (with the following axioms) completely axiomatizes both of the structures $\langle \mathbb{R}; < \rangle$ and $\langle \mathbb{Q}; < \rangle$.*

(**A**$_<$) $\forall x,y\,(x<y \rightarrow y \not< x)$
(**T**$_<$) $\forall x,y\,(x<y<z \rightarrow x<z)$
(**L**$_<$) $\forall x,y\,(x<y \vee x=y \vee y<x)$
(**D**$_<$) $\forall x,y\,(x<y \rightarrow \exists w\,[x<w<y])$
(**U**$_<$) $\forall x \exists u\,(x<u)$
(**B**$_<$) $\forall x \exists v\,(v<x)$

For a proof, note that the axioms are true in $\langle \mathbb{R}; < \rangle$ and $\langle \mathbb{Q}; < \rangle$; so, it suffices to show that the above theory admits QE. For that we use Lemma 4 and show the equivalence of every formula of the form $\exists x \bigwedge_i \gamma_i(x)$ to a quantifier-free formula, where each $\gamma_i$ is an atom or negated atom. The equivalences $\neg(a<b) \leftrightarrow (a=b) \vee (b<a)$ and $\neg(a=b) \leftrightarrow (a<b) \vee (b<a)$, which are provable in the theory, allow us to neglect negated atomic formulas. Thus, we need to eliminate the quantifier of the formulas of the form $\exists x(\bigwedge_i u_i < x \wedge \bigwedge_j x < v_j \wedge \bigwedge_k x = w_k)$ only—note that $x=x$ is equivalent to $\top$, and $x<x$ to $\bot$. But that formula is equivalent to $\bigwedge_i u_i < w_0 \wedge \bigwedge_j w_0 < v_j \wedge \bigwedge_k w_0 = w_k$, if the conjunction $\bigwedge_k x = w_k$ is non-empty, and to $\bigwedge_{i,j} u_i < v_j$, if it is empty (non-existent) and none of the other conjunctions are empty; if any of $\bigwedge_i u_i < x$ or $\bigwedge_j x < v_j$ is also empty, then the original formula is equivalent to $\top$.

The order relation behaves very differently on $\mathbb{Z}$ and $\mathbb{N}$, since here it is a *discrete* order, in the sense that every element has an immediate successor. Let us denote the successor function $x \mapsto (x+1)$ by $\mathfrak{s}$; and let $(x \leqslant y)$ abbreviate $(x<y) \vee (x=y)$. For a proof of the following theorem, first proved by A. Robinson and E. Zakon (1960, Theorem 2.12), see e.g. [1, Theorem 2].

**THEOREM 6 (An Axiomatization for $\langle \mathbb{Z}; < \rangle$)** *The (finitely axiomatized) theory of discrete linear orders without endpoints completely axiomatizes the structure $\langle \mathbb{Z}; <, \mathfrak{s} \rangle$; this theory consists of the axioms **A**$_<$, **T**$_<$, **L**$_<$ (Theorem 5) along with*

(**S**$_<$) $\forall x,y\,(x<y \leftrightarrow \mathfrak{s}(x) \leqslant y)$
(**P**$_<$) $\forall x \exists w\,(\mathfrak{s}(w)=x)$

The following has been proved in e.g. [6, Theorem 32A].

**THEOREM 7 (An Axiomatization for $\langle \mathbb{N}; < \rangle$)** *The (finitely axiomatizable) theory of discrete linear orders with the least element and without the last element completely axiomatizes the structure $\langle \mathbb{N}; \mathbf{0}, \mathfrak{s}, < \rangle$; this theory consists of the axioms **A**$_<$, **T**$_<$, **L**$_<$ (Theorem 5) together with **S**$_<$ (Theorem 6) and*

(**Z**$_<$) $\forall x\,(\mathbf{0} \leqslant x)$
(**P**$_<^0$) $\forall x \exists w\,(\mathbf{0} < x \rightarrow \mathfrak{s}(w)=x)$

Let us note that $\forall x\,[x < \mathfrak{s}(x)]$ is provable from **S**$_<$; and so one can show that $\forall x,y\,[x<y \leftrightarrow \mathfrak{s}(x) < \mathfrak{s}(y)]$ follows from **S**$_<$, **T**$_<$ and **L**$_<$. Therefore, Peano's axioms $\forall x\,(\mathfrak{s}(x) \neq \mathbf{0})$ and $\forall x,y\,(\mathfrak{s}(x)=\mathfrak{s}(y) \rightarrow x=y)$ are provable from the axiom system **A**$_<$, **T**$_<$, **L**$_<$, **Z**$_<$, and **S**$_<$.

We now study the addition operation ($+$) in number systems. The most obvious properties of addition are associativity and commutativity (see **A**$_+$ and **C**$_+$ in Theorem 8). Of course, in all of our number systems there is an additive unit element (zero $\mathbf{0}$), and in all but one (the natural numbers) every element has an additive inverse (the minus element). In $\mathbb{C}$, $\mathbb{R}$, and $\mathbb{Q}$ addition is torsion-free and divisible (see **T**$_+$ and **D**$_+$ in Theorem 8); it is hard to find any other property of $+$ in $\mathbb{C}, \mathbb{R}, \mathbb{Q}$ that does not follow from the above-mentioned properties.

For axiomatizing the structures $\langle \mathbb{C}; + \rangle$, $\langle \mathbb{R}; + \rangle$, and $\langle \mathbb{Q}; + \rangle$ we add the constant symbol $\mathbf{0}$ and the unary function symbol $-$ to the language; needless to say, $n \cdot x$ abbreviates the expression $x+\cdots+x$ ($n$ times) for $n \in \mathbb{N}$.

**THEOREM 8 (An Axiomatization for $\langle\mathbb{C};+\rangle$, $\langle\mathbb{R};+\rangle$, $\langle\mathbb{Q};+\rangle$)**
*The first-order theory of non-trivial, divisible, torsion-free, and commutative groups (with the following infinite set of axioms) completely axiomatizes the structures $\langle\mathbb{Q};0,-,+\rangle$, $\langle\mathbb{R};0,-,+\rangle$, and $\langle\mathbb{C};0,-,+\rangle$.*

- **($A_+$)** $\forall x,y,z\,(x+(y+z)=(x+y)+z)$
- **($C_+$)** $\forall x,y\,(x+y=y+x)$
- **($U_+$)** $\forall x\,(x+\mathbf{0}=x)$
- **($I_+$)** $\forall x\,(x+(-x)=\mathbf{0})$
- **($N_+$)** $\exists u\,(u\neq\mathbf{0})$
- **($T_+$)** $\{\forall x\,(n\,.\,x=\mathbf{0}\to x=\mathbf{0})\}_{n>0}$
- **($D_+$)** $\{\forall x\exists v\,(x=n\,.\,v)\}_{n>0}$

We show that the theory admits QE by using Lemma 4. Every atomic formula in the language $\{\mathbf{0},-,+\}$ that contains $x$ can be equivalently written in the form $n\,.\,x=t$ for some $n\in\mathbb{N}^+$ and some $x$–free term $t$. By $a=b\leftrightarrow k\,.\,a=k\,.\,b$, which is provable from the above axioms, it suffices to eliminate the quantifier of $\exists x(\bigwedge_i q\,.\,x=t_i\wedge\bigwedge_j q\,.\,x\neq s_j)$, which by $D_+$ (for $n=q$) is equivalent to $\exists y(\bigwedge_i y=t_i\wedge\bigwedge_j y\neq s_j)$. Now, if the conjunct $\bigwedge_i y=t_i$ is nonempty, then this is equivalent to $\bigwedge_i t_0=t_i\wedge\bigwedge_j t_0\neq s_j$, and if $\bigwedge_i y=t_i$ is empty, then it is equivalent to $\top$, since by $N_+$ there are infinitely many members (for any $u\neq\mathbf{0}$ we have $n\,.\,u\neq m\,.\,u$ for every $n\neq m$).

The axiomatization of $\langle\mathbb{Z};+\rangle$ illustrates a case that one might need to substantially enrich the language of the structure to have QE. As an example, the formula $\exists v\,(x=v+v)$, stating that $x$ is even, is not equivalent to any quantifier-free formula in $\langle\mathbb{Z};0,-,+\rangle$. However, if we add the binary relation symbol $\equiv_2$ of congruence modulo 2 to the language, then that formula will be equivalent to $x\equiv_2\mathbf{0}$.

The quantifier elimination of the theory of the structure $\langle\mathbb{Z};\mathbf{0},\mathbf{1},\{\equiv_n\}_{n>1},-,+\rangle$ can be shown by using a generalized form of the Chinese Remainder Theorem in Number Theory. The Chinese remainder theorem says that a given system of congruence equations $\{x\equiv_{n_i}r_i\}_{i<N}$ has a solution (in $\mathbb{Z}$) if $n_i$ and $n_j$ are coprime for every $i<j<N$. The generalized Chinese remainder theorem says that the system $\{x\equiv_{n_i}r_i\}_{i<N}$ of congruence equations has a solution if and only if for every $i<j<N$ we have $r_i\equiv_{d_{i,j}}r_j$, where $d_{i,j}$ is the greatest common divisor of $n_i$ and $n_j$. Since such systems either have no solution or have infinitely many solutions, then we can state this more general theorem as follows.

**PROPOSITION 9 (General Chinese Remainder Theorem)**
*If $n_i>1$ for every $i<N$, then for every $\{r_i\}_{i<N}$ and $\{s_j\}_{j<M}$,*

$$\exists x(\bigwedge_{i<N}x\equiv_{n_i}r_i\wedge\bigwedge_{j<M}x\neq s_j)\iff\bigwedge_{i<N}r_i\equiv_{d_{i,j}}r_j,$$

*where $d_{i,j}$ is the greatest common divisor of $n_i$ and $n_j$.*

For three different proofs of Proposition 9, which is a kind of QE by itself, see [20, Propositions 4.5 and 4.1] and [1, Proposition 2] which are due to Ore (1951), Mahler (1958) and Fraenkel (1963) respectively.

We add the congruence relations $\equiv_n$ modulo every natural $n>1$, along with the constant $\mathbf{1}$, to the language; let $\bar{i}$ abbreviate $\mathbf{1}+\cdots+\mathbf{1}$ ($i$ times) for every $i\in\mathbb{N}$.

**THEOREM 10 (An Axiomatization for $\langle\mathbb{Z};+\rangle$)** *The theory whose axioms are $A_+$, $C_+$, $U_+$, $I_+$, and $T_+$ (Theorem 8) with the following axioms completely axiomatizes the structure $\langle\mathbb{Z};\mathbf{0},\mathbf{1},\{\equiv_n\}_{n>1},-,+\rangle$.*

- **($E^+$)** $\{\forall x,y\,[x\equiv_n y\leftrightarrow\exists u(x=y+n\,.\,u)]\}_{n>1}$
- **($E_+$)** $\{\forall x\,[\bigvee_{i<n}(x\equiv_n\bar{i})]\}_{n>1}$
- **($E'_+$)** $\{\bigwedge_{0<i<n}(\bar{i}\not\equiv_n\mathbf{0})\}_{n>1}$

For showing that the theory admits QE by Lemma 4, we note that every atomic formula of $x$ in $\{\mathbf{0},\mathbf{1},-,+\}\cup\{\equiv_n|\,n>1\}$ is equivalent to either $m\,.\,x=t$ or $m\,.\,x\equiv_n t$ for some $m,n\in\mathbb{N}^+$ and some $x$–free term $t$. By the provable equivalence $(a\not\equiv_n b)\leftrightarrow\bigvee_{0<i<n}(a\equiv_n b+\bar{i})$ it suffices to show that the formula

$$\exists x(\bigwedge_i q_i\,.\,x\equiv_{n_i}r_i\wedge\bigwedge_j q_j\,.\,x\neq s_j\wedge\bigwedge_k q_k\,.\,x=t_k)$$

is equivalent to a quantifier-free formula. From the provable equivalences $(a=b)\leftrightarrow(k\,.\,a=k\,.\,b)$ and $(a\equiv_n b)\leftrightarrow(k\,.\,a\equiv_{kn}k\,.\,b)$ we can assume that all the $q_i$'s, $q_j$'s and $q_k$'s are equal, to say $q$. Then, the above formula is equivalent to

$$\exists y(y\equiv_q\mathbf{0}\wedge\bigwedge_i y\equiv_{n_i}r_i\wedge\bigwedge_j y\neq s_j\wedge\bigwedge_k y=t_k).$$

We can assume that the conjunct $\bigwedge_k y=t_k$ is empty (see the proofs of Theorems 5,8); now the result immediately follows from Proposition 9 (which is provable from the axioms).

As for $\mathbb{N}$, even the language $\{\mathbf{0},\mathbf{1},-,+\}\cup\{\equiv_n|\,n>1\}$ is not sufficiently rich for QE, as the formula $\exists v(x+v=y)$ is not equivalent to a quantifier-free formula (it is equivalent to $x\leqslant y$) in $\mathbb{N}$. Here, QE is possible when we add the order relation to the language.

**THEOREM 11 (An Axiomatization for $\langle\mathbb{N};<,+\rangle$)** *The theory with the axioms $A_<$, $T_<$, $L_<$ (Theorem 5), $S_<$ (Theorem 6), $Z_<$, $P_<^0$ (Theorem 7), $A_+$, $C_+$, $U_+$ (Theorem 8), $E^+$, $E_+$ (Theorem 10) with the following axioms completely axiomatizes the structure $\langle\mathbb{N};\mathbf{0},\mathbf{1},<,\{\equiv_n\}_{n>1},+\rangle$.*

- **($M_+$)** $\forall x,y\,(x<y\to\exists v\,[x+v=y])$
- **($O_+$)** $\forall x,y,z\,(x<y\to x+z<y+z)$

A proof of Theorem 11 can be found in [6, Theorem 32E] (without presenting an explicit axiomatization; though one can see that the proof goes through with our suggested axioms). For a proof of the following theorem see e.g. [1, Theorem 5]; other proofs can be found in [13, § 4.III] and [23, §§ III.4.2].

**THEOREM 12 (An Axiomatization for $\langle\mathbb{Z};<,+\rangle$)** *The theory with the axioms $A_<$, $T_<$, $L_<$ (Theorem 5), $S_<$, $P_<$ (Theorem 6), $A_+$, $C_+$, $U_+$, $I_+$ (Theorem 8), $E_+$ (Theorem 10) and finally $O_+$ (Theorem 11), with $\mathfrak{s}(x)$ set to $x+\mathbf{1}$, completely axiomatizes the structure $\langle\mathbb{Z};\mathbf{0},\mathbf{1},<,\{\equiv_n\}_{n>1},-,+\rangle$.*

The following theorem (stating that the order and addition structure of rational and real numbers can be axiomatized by the theory of non-trivial divisible commutative ordered groups) can be proved by combining the techniques of the proofs of Theorems 5 and 8 (cf. [1, Theorem 4]).

**THEOREM 13 (Axiomatizing $\langle \mathbb{Q}; <, + \rangle$ and $\langle \mathbb{R}; <, + \rangle$)** *The theory with the axioms $A_<$, $T_<$, $L_<$ (Theorem 5), $A_+$, $C_+$, $U_+$, $I_+$, $N_+$, $D_+$ (Theorem 8), and finally $O_+$ (Theorem 11) completely axiomatizes the structures $\langle \mathbb{Q}; 0, <, -, + \rangle$ and $\langle \mathbb{R}; 0, <, -, + \rangle$.*

Let us note that the axioms $D_<$, $U_<$, $B_<$ (in Theorem 5) and $T_+$ (in Theorem 8) are provable from the axiom system presented in Theorem 13 just the way that are proved in classical analysis.

## § Number Systems (Addition & Multiplication).

**DEFINITION 14 (Field)** A *field* is a structure over $\{0, 1, +, -, \times, ^{-1}\}$ that satisfies $A_+$, $C_+$, $U_+$, $I_+$ (Theorem 8), and the following axioms:

- $(A_\times)$ $\forall x, y, z \, (x \cdot (y \cdot z) = (x \cdot y) \cdot z)$
- $(C_\times)$ $\forall x, y \, (x \cdot y = y \cdot x)$
- $(U_\times)$ $\forall x \, (x \cdot 1 = x)$
- $(I_\times)$ $\forall x \, (x \neq 0 \rightarrow x \cdot x^{-1} = 1)$
- $(D_\times)$ $\forall x, y, z \, [x \cdot (y+z) = (x \cdot y) + (x \cdot z)]$

A field has *characteristic zero* if it moreover satisfies

- $(C_0)$ $\{\overline{n} \neq 0\}_{n>0}$

where, as we recall, $\overline{n}$ abbreviate $1 + \cdots + 1$ ($n$ times). $\diamond$

The field $\langle \mathbb{C}; +, \times \rangle$ is well known to be algebraically closed since it satisfies the Fundamental Theorem of Algebra, i.e., it has a root for every non-trivial polynomial (with coefficients in $\mathbb{C}$). It can be even said that it was created for having all the roots of the polynomials (with real or complex coefficients). This is all one can say about the complex field in the first-order setting, since the theory of algebraically closed fields of characteristic zero is complete, and so it axiomatizes $\langle \mathbb{C}; +, \times \rangle$; see also [15]. The following result was proved by Tarski (1936); see e.g. [13, § 4.IV] for a proof.

**THEOREM 15 (An Axiomatization for $\langle \mathbb{C}; +, \times \rangle$)** *The theory of algebraically closed fields of characteristic zero, with the axioms $A_+$, $C_+$, $U_+$, $I_+$, $A_\times$, $C_\times$, $U_\times$, $I_\times$, $D_\times$, $C_0$ (Definition 14) along with the following axioms, completely axiomatizes the structure $\langle \mathbb{C}; 0, 1, -, +, \times, ^{-1} \rangle$.*

- $(\mathit{FTA}_\mathbb{C})$ $\{\forall \langle a_i \rangle_{i<n} \exists x \, (x^n + \sum_{i<n} a_i x^i = 0)\}_{n>1}$

For super-careful readers let us note that (i) every non-trivial polynomial can be taken to be a monic by dividing it with the leading (non-zero) coefficient; (ii) the multiplicative inversion ($x \mapsto x^{-1}$) is not really a total function, since it is not defined on zero, but one can make the convention $0^{-1} = 0$ without any danger; (iii) and finally, $x^i$ abbreviates the algebraic expression $x \times \cdots \times x$ ($i$ times) of course.

For studying the structure $\langle \mathbb{R}; +, \times \rangle$ we first note that order is definable in it: $u \leqslant v \iff \exists x (u + x^2 = v)$; and $\langle \mathbb{R}; <, +, \times \rangle$ is an ordered field, see e.g. [17]. An ordered field satisfies the order axioms $A_<$, $T_<$, $L_<$ (Theorem 5), the axioms of fields (in Definition 14), $O_+$ (Theorem 11), and $O_\times$ (Theorem 16 below).

Of course this is not all one can say about $\langle \mathbb{R}; <, +, \times \rangle$. On the other hand, not much can one say about it in the first-order framework; only that every positive real number has a square root, and every polynomial of even degree can be factorized into some quadratic polynomials (see $\mathit{FTA}_\mathbb{R}$ in Theorem 16). This last statement is indeed equivalent to (a real version of) the fundamental theorem of algebra. As some examples, let us note quadratic factorizations of the following quartics:

$$x^4 + 1 = (x^2 + \sqrt{2}x + 1)(x^2 - \sqrt{2}x + 1),$$
$$x^4 - x^2 + 1 = (x^2 + \sqrt{3}x + 1)(x^2 - \sqrt{3}x + 1), \text{ and}$$
$$x^4 - x + 1 = (x^2 + \sqrt{r}x + \tfrac{2\sqrt{r}}{r\sqrt{r}-1})(x^2 - \sqrt{r}x + \tfrac{r\sqrt{r}-1}{2\sqrt{r}}),$$

where $r$ is the unique positive real number that satisfies the cubic equation $r^3 - 4r - 1 = 0$. For a proof of the following result of Tarski (1936) see e.g. [21, Appendix] which is a modified version of the proof presented in [13, § 4.V].

**THEOREM 16 (An Axiomatization for $\langle \mathbb{R}; <, +, \times \rangle$)** *Theory of real closed ordered fields which is axiomatized by $A_<$, $T_<$, $L_<$ (Theorem 5), the axioms of fields (Definition 14), and $O_+$ (Theorem 11) along with the following axioms completely axiomatizes the structure $\langle \mathbb{R}; 0, 1, <, -, +, \times, ^{-1} \rangle$.*

- $(O_\times)$ $\forall x, y, z \, (0 < z \wedge x < y \rightarrow x \cdot z < y \cdot z)$
- $(S_\times)$ $\forall x \, (0 < x \rightarrow \exists u \, [x = u^2])$
- $(\mathit{FTA}_\mathbb{R})$ $\{\forall \langle a_i \rangle_{i<2n} \exists \langle b_j, c_j \rangle_{j<n} \forall x$
  $[(x^{2n} + \sum_{i<2n} a_i x^i) = \prod_{j<n}(x^2 + b_j x + c_j)]\}_{n>1}$

We note that by $S_\times$ (and the axioms of ordered fields) the high-school equivalence for the existence of the roots of quadratic polynomials can be proved:

$$\exists x (x^2 + bx + c = 0) \leftrightarrow \exists x [(2x+b)^2 = b^2 - 4c] \leftrightarrow b^2 \geqslant 4c.$$

It can be easily seen that $\mathit{FTA}_\mathbb{C}$ (in Theorem 15) is equivalent to the statement that every monic is equal to a product of some linear polynomials:

$\mathit{FTA}_\mathbb{C} \equiv \{\forall \langle a_i \rangle_{i<n} \exists \langle b_j \rangle_{j<n} \forall x$
$[(x^n + \sum_{i<n} a_i x^i) = \prod_{j<n}(x + b_j)]\}_{n>1}$

which resembles $\mathit{FTA}_\mathbb{R}$ (in Theorem 16). Let us note a couple of consequences of $\mathit{FTA}_\mathbb{R}$ (from [21]):

**PROPOSITION 17 ($\mathit{FTA}_\mathbb{R} \implies \mathit{RCF} + \mathit{IVT}$)** *If every even-degree polynomial can be factorized into some quadratic polynomials in an ordered field in which every positive element has a square root, then every odd-degree polynomial has a root and the polynomial intermediate value theorem holds in it.*

For a proof suppose that the polynomial $\mathfrak{p}(x)$ is of degree $m$ and $\mathfrak{p}(u)\mathfrak{p}(v) < 0$ holds for some $u < v$. Put $\mathfrak{q}(x)$ to be the polynomial $\tfrac{1}{\mathfrak{p}(u)}(1+x^2)^m \mathfrak{p}(u + \tfrac{v-u}{1+x^2})$; then $\mathfrak{q}(x) = x^{2m} + \mathfrak{r}(x^2)$ for some polynomial $\mathfrak{r}(x)$ with degree less than $m$. So, $\mathfrak{q}(x)$ can be factorized to say $\prod_{j<m}(x^2 + b_j x + c_j)$. Now we have $\prod_{j<m} c_j = \mathfrak{q}(0) = \tfrac{\mathfrak{p}(v)}{\mathfrak{p}(u)} < 0$ and so $c_j < 0$ for some $j$; then we have $b_j^2 > 4c_j$ and so the quadratic $x^2 + b_j x + c_j = 0$ has a root, such as $s$. Now, $r = u + \tfrac{v-u}{1+s^2}$ is a root of $\mathfrak{p}(x) = 0$ that satisfies $u < r < v$. By a classical real analytic argument, if the intermediate value theorem holds for polynomials in an ordered field, then every odd-degree polynomial has a root in it.

So, the fundamental theorem of algebra is really *fundamental* since it can prove some basic theorems in algebra, and it is a kind of fundamental theorem for the mathematical analysis of polynomials as well; see [21] for more details.

So far, we have observed two applications of mathematics (especially number theory and algebra) in mathematical logic:

1. The (Generalized) Chinese Remainder Theorem
2. The Fundamental Theorem of Algebra

1. Proposition 9 was used in proving that the axiomatic system suggested for the additive structure of integer numbers $\langle\mathbb{Z};+\rangle$ has QE and so it is a complete theory (Theorem 10). It is worth noting that Gödel (1931, Lemma 1) also had used the (non-generalized) Chinese remainder theorem in his proof of the first incompleteness theorem for the coding technicalities.

2. The truly *fundamental theorem* of elementary algebra and elementary analysis was used for axiomatizing the additive and multiplicative structures of complex and real numbers, $\langle\mathbb{C};+,\times\rangle$ and $\langle\mathbb{R};+,\times\rangle$, noting that order ($<$) is definable in $\langle\mathbb{R};+,\times\rangle$ (Theorems 15 and 16).

Now, we present two applications of mathematical logic in other areas of mathematics (especially algebraic geometry):

I. The Tarski-Seidenberg Principle
II. Hilbert's 17th Problem

I. Theorem 16, like many other theorem of QE, is proved by using Lemma 4. Let us see how the proof can proceed: first, we note that all the atomic formulas of $x$ over the language $\{\mathbf{0},\mathbf{1},-,+,\times,{}^{-1},<\}$ are equivalent to $p(x)=0$ or $p(x)>0$ for a polynomial $p$. Second, negation can be eliminated (see the proof of Theorem 5), so QE over this language is equivalent to "the existence of a solution of a system of polynomial equations and inequalities is equivalent to a system of some equations and inequalities between the coefficients of those polynomials". As an example, $\exists x(ax^2+bx+c=0)$ is known to be equivalent to $(a^2>0 \wedge b^2 \geqslant 4ac) \vee (a=0 \wedge b^2>0) \vee (a=0 \wedge b=0 \wedge c=0)$. This is called the Tarski-Seidenberg principle in real algebraic geometry (see [3, §1.4]), which is exactly what the translation of Lemma 4 would be in the proof of Theorem 16 (cf. [15]).

II. Hilbert's celebrated 17th Problem asked (see e.g. [22]): *Given a multivariate polynomial that takes only non-negative values over the reals, can it be represented as a sum of squares of rational functions?* Let us note a couple of examples:

$$x^4-x+1 = (x^2-\tfrac{1}{2})^2+(x-\tfrac{1}{2})^2+(\tfrac{1}{\sqrt{2}})^2 \text{ and}$$

$$(x^2+y^2)^2[x^4y^2+x^2y^4+1-3x^2y^2] = (x^2-y^2)^2 + [x^2y(x^2+y^2-2)]^2+[xy^2(x^2+y^2-2)]^2+[xy(x^2+y^2-2)]^2.$$

A consequence of the Tarski-Seidenberg principle is the Artin-Lang Homomorphism theorem [3, Theorem 4.1.2] which gives a positive answer to the problem; see [3, Theorem 6.1.1]. Let us note that by the fundamental theorem of algebra every non-negative polynomial of one variable can be written as a sum of the squares of some polynomials;[1] but there are non-negative polynomials of two variables that cannot be written as such. One example (see [22]) is Motzkin (1969)'s polynomial $x^4y^2+x^2y^4+1-3x^2y^2$; of course it is the sum of the squares of some *rational* functions (see the second example above).

The next structures that we study over the language $\{+,\times\}$ are $\mathbb{Q}$, $\mathbb{Z}$, and $\mathbb{N}$. Here the story becomes dramatically different. To start with, let us note that the axiomatic systems presented for the ordered structures $\langle\mathbb{N};<\rangle$, $\langle\mathbb{Z};<\rangle$, $\langle\mathbb{Q};<\rangle$, and $\langle\mathbb{R};<\rangle$ were all finite (Theorems 5,6,7). Other axiomatic systems were not finite, but were presented in a way that one can recognize whether a given sentence is an axiom of that system or not, in the sense that a properly designed algorithm can recognize them. In the other words, the axiomatic theories for the studied structures were decidable by an algorithm.

---
[1] The sum of squares for a given polynomial may not be unique, as the identity $(4x^2+1)^2 = (4x)^2+(4x^2-1)^2$ shows.

To make precise the forthcoming definition, let us make the convention that all our first-order individual variables are $\vartheta, \vartheta', \vartheta'', \vartheta''', \cdots$, made up from $\vartheta$ and $'$. Let us fix the following finite set of symbols as an alphabet: $\mathcal{A} = \{\neg, \wedge, \vee, \forall, \exists, (,), \vartheta, ', \mathbf{0}, \mathbf{1}, <, =, +, -, \times, \mathfrak{exp}\}$. Every formula over the first-order language $\{\mathbf{0}, \mathbf{1}, <, +, -, \times, \mathfrak{exp}\}$ is a string (i.e., a finite sequence) of the elements of $\mathcal{A}$. There exists an algorithm that decides (outputs yes or no) if a given such string as input is a well-founded formula or not.

**DEFINITION 18 (Decidability)** A set $B$ of strings of symbols from $\mathcal{A}$ is *decidable* when there exists an algorithm such that for a given string as input outputs yes if it belongs to $B$ and outputs no otherwise. ⋄

Let us note that we have not fixed a rigorous definition for the informal notion of *algorithm* in the above definition; it could be a *recursive function* or a *Turing machine* (see [2]). By the Church-Turing thesis all such formally rigorous and equivalent definitions do define the informal notion of algorithm; so we do not need to fix a formalization. "Axiomatizable" usually means *axiomatizable by a decidable set of axioms*; though more often the decidability of the axiom set is not explicitly mentioned.

**DEFINITION 19 (Axiomatizability)** A theory or a structure is said to be *axiomatizable* when there exists a decidable set of sentences that completely axiomatizes it. ⋄

All the theories and structures that we have studied so far are axiomatizable by a decidable set of sentences. Actually, a structure is axiomatizable by a decidable set of sentences if and only if it has a decidable theory; see e.g. [6, Corollary 26I].

Decidability implies axiomatizability, since one only needs to algorithmically list all the sentences and pick the ones that hold true; thus a decidable set of axioms is obtained. If $\mathfrak{A}$ is axiomatizable, then for a given sentence $\psi$ run this algorithm for consecutive $n$'s starting from $n=1$: list all the theorems that are proved in $n$ steps or less from the first $n$ axioms (if $n$ exceeds the number of axioms, then use all the finitely many axioms); if $\psi$ or $\neg\psi$ appears in the list, then output yes or no accordingly. The algorithm will surely terminate (for some $n$) since the axiomatic system *completely* axiomatizes $\mathfrak{A}$.

Now, the shocking result of Gödel's incompleteness theorem (1931) is that the structure $\langle\mathbb{N};+,\times\rangle$ is not axiomatizable. As the history goes, Presburger (1929) proved the axiomatizability of $\langle\mathbb{N};+\rangle$ and Skolem (1930) announced the axiomatizability of $\langle\mathbb{N};\times\rangle$ (see [23]); so $\langle\mathbb{N};+,\times\rangle$ was expected to be axiomatizable, that would confirm Hilbert's Programme (see e.g. [5, 7]).

**THEOREM 20 (Non-Axiomatizability of $\langle\mathbb{N};+,\times\rangle$)** *The full first-order theory of $\langle\mathbb{N};+,\times\rangle$ is not axiomatizable by any decidable set of sentences.*

Of course, there does exist an undecidable set of sentences that completely axiomatizes $\langle\mathbb{N};+,\times\rangle$; that is the so-called *true* *arithmetic*, the set of all the sentences that are true in $\mathbb{N}$. The non-axiomatizability of $\langle\mathbb{Z};+,\times\rangle$ is inherited from $\langle\mathbb{N};+,\times\rangle$ since the set $\mathbb{N}$ is definable in $\langle\mathbb{Z};+,\times\rangle$ by Lagrange's four square theorem (see e.g. [23, Theorem II.3.8]): let $\mathcal{N}(x)$ be the formula $\exists u,v,w,z(u^2+v^2+w^2+z^2 = x)$. Then for every $m \in \mathbb{Z}$ we have [$m \in \mathbb{N}$ if and only if $\mathcal{N}(m)$ is true in $\mathbb{Z}$].

For every formula $\varphi$ over $\{+,\times\}$, let $\varphi^\mathcal{N}$ result from $\varphi$ by changing every $\forall x\,\Theta$ to $\forall x\,[\mathcal{N}(x) \to \Theta]$ and $\exists x\,\Theta$ to $\exists x\,[\mathcal{N}(x) \wedge \Theta]$; that is *relativizing* all the bounded variables to $\mathcal{N}$. Now, for every sentence $\theta$ over $\{+,\times\}$ we have: $\theta$ is true in $\langle\mathbb{N};+,\times\rangle$ if and only if $\theta^\mathcal{N}$ is true in $\langle\mathbb{Z};+,\times\rangle$.

So, it follows that the structure $\langle\mathbb{Z};+,\times\rangle$ is not axiomatizable by any decidable set of sentences $T$, since otherwise $\langle\mathbb{N};+,\times\rangle$ would be axiomatizable by the decidable set of sentences $T^\mathcal{N} = \{\theta^\mathcal{N} \mid \theta \in T\}$. For another definition of $\mathbb{N}$ in $\langle\mathbb{Z};+,\times\rangle$ see Robinson's paper [19] where it is proved that $\mathbb{Z}$ is definable in $\langle\mathbb{Q};+,\times\rangle$ as well (see also [18]).

**COROLLARY 21 (On $\langle\mathbb{Z};+,\times\rangle$ and $\langle\mathbb{Q};+,\times\rangle$)** *The structures $\langle\mathbb{Z};+,\times\rangle$ and $\langle\mathbb{Q};+,\times\rangle$ are not axiomatizable.*

§ **NUMBER SYSTEMS (Multiplication & Exponentiation).**

We saw that by Tarski's result $\langle\mathbb{C};+,\times\rangle$ is axiomatizable; then its theory is decidable, and so is the theory of $\langle\mathbb{C};\times\rangle$. Thus, $\langle\mathbb{C};\times\rangle$ is axiomatizable by a decidable set of sentences; but what is that axiomatic system? This question was answered in [20, Theorem 2.2] by providing an explicit axiomatization for the multiplicative structure of complex numbers:

**THEOREM 22 (An Axiomatization for $\langle\mathbb{C};\times\rangle$)** *The structure $\langle\mathbb{C};0,1,\{\omega_n\}_{n>1},\times,{}^{-1}\rangle$ is axiomatizable by $A_\times$, $C_\times$, $U_\times$, $I_\times$ (Theorem 16) along with the following axioms:*

$(Z_\times)$ $\forall x\,(x \cdot 0 = 0 = 0^{-1})$
$(D_\times)$ $\{\forall x\exists v\,(x = v^n)\}_{n>0}$
$(R_\times)$ $\{\forall x\,[x^n = 1 \leftrightarrow \bigvee_{i<n} x = (\omega_n)^i]\}_{n>1}$
$(R^\times)$ $\{\bigwedge_{i<j<n}(\omega_n)^i \neq (\omega_n)^j\}_{n>1}$

*where $\omega_n$ is interpreted as $[\cos(2\pi/n)+i\sin(2\pi/n)]$ for every $n>1$; thus $\omega_2 = -1$, $\omega_3 = (-1/2) + i(\sqrt{3}/2)$ and $\omega_4 = i$.*

The same question can be asked about the real numbers: we know that $\langle\mathbb{R};\times\rangle$ is decidable by Tarski's result that $\langle\mathbb{R};+,\times\rangle$ is axiomatizable; but what is an explicit axiomatization for $\langle\mathbb{R};\times\rangle$? For its answer we need to add the *positivity predicate*, denoted $\mathcal{P}(x)$, to the language. The following result is proved in [20, Theorem 3.3].

**THEOREM 23 (An Axiomatization for $\langle\mathbb{R};\times\rangle$)** *The structure $\langle\mathbb{R};0,1,-1,\mathcal{P},\times,{}^{-1}\rangle$ is axiomatizable by $A_\times$, $C_\times$, $U_\times$, $I_\times$, $Z_\times$ (Theorem 22) along with the following axioms:*

$(N_\times)$ $\exists u\,(u \neq 0,1,-1)$
$(D_\times^o)$ $\{\forall x\exists v\,(x = v^{2n+1})\}_{n>0}$
$(R_\times^e)$ $\{\forall x\,(x^{2n} = 1 \leftrightarrow x = 1 \vee x = -1)\}_{n>1}$
$(P)$ $\forall x\,(\mathcal{P}(x) \leftrightarrow \exists y \neq 0[x = y^2])$
$(P_\times)$ $\forall x, y \neq 0(\mathcal{P}(xy) \longleftrightarrow [\mathcal{P}(x) \leftrightarrow \mathcal{P}(y)])$
$(P_\times^-)$ $\forall x \neq 0[\neg\mathcal{P}(x) \leftrightarrow \mathcal{P}([-1]x)]$

Let us note that the multiplicative structure of positive real numbers $\langle\mathbb{R}^+;\times\rangle$ is a non-trivial, divisible, torsion-free, and commutative group, since it is isomorphic to $\langle\mathbb{R};+\rangle$ via the mapping $x \mapsto \ln(x)$.

For axiomatizing the multiplicative structure of rational numbers $\langle\mathbb{Q};\times\rangle$ we first axiomatize the multiplicative structure of positive rational numbers $\langle\mathbb{Q}^+;\times\rangle$ noting that one can obtain an axiomatization for $\langle\mathbb{Q};\times\rangle$ by adding the constants $0,-1$ and the predicate $\mathcal{P}(x)$ to the language and adding $Z_\times$, $N_\times$, $P$, $P_\times$, and $P_\times^-$ (Theorem 23) to the axioms. The following is proved in [20, Theorem 4.11]:

**THEOREM 24 (An Axiomatization for $\langle\mathbb{Q}^+;\times\rangle$)** *The first-order structure $\langle\mathbb{Q}^+;1,\times,{}^{-1}\rangle$ is axiomatizable by $A_\times$, $C_\times$, $U_\times$, $I_\times$ (Definition 14) along with the following axioms:*

$(T_\times)$ $\{\forall x\,(x^n = 1 \to x = 1)\}_{n>1}$
$(M_\times)$ $\{\forall\langle x_i\rangle_{i<k}\exists v \forall y \bigwedge_{i<k}(v^n x_i \neq y^{m_i})\}_{n,k}$
*where $n, k \in \mathbb{N}$, and no $m_i \in \mathbb{N}$ divides $n$.*

The axioms $M_\times$ in Theorem 24 state that for every sequence $x_0, \ldots, x_{k-1}$ of positive rational numbers and every sequence $m_0, \ldots, m_{k-1}$ of natural numbers none of which divides the natural number $n$, there exists a positive rational number $v$ such that for every $i<k$ none of $v^n x_i$'s is an $m_i$-power of a positive rational number. To see that this holds in $\mathbb{Q}^+$ it suffices to take $v$ to be a prime number that does not divide the numerators and denominators of any of $x_i$'s. This does not hold if some $m_i$ divides $n$ since $x_i$ could be an $m_i$th power; it does not hold in $\mathbb{R}^+$ either since every positive real number has an $m_i$th root.

The next structures that we study over the language $\{\times\}$ are $\mathbb{Z}$ and $\mathbb{N}$. Here too, as we saw, it suffices to study $\langle\mathbb{N}^+;\times\rangle$ first, and then for $\langle\mathbb{N};\times\rangle$ we need to add $0$ and the axiom $Z_\times$, and for $\langle\mathbb{Z};\times\rangle$ we need to add $-1,\mathcal{P}$ and the axioms $P_\times$ and $P_\times^-$. Since studying the axioms of $\langle\mathbb{N}^+;\times\rangle$ will not be needed later, and they are too many to be listed in the main body of the paper, and explaining them will take much time and will distract the flow of the paper, we apologetically postpone it to the Appendix.

Let us move on to the language $\{<,\times\}$ over which $\mathbb{R}$ and $\mathbb{Q}$ are axiomatizable, while $\mathbb{Z}$ and $\mathbb{N}$ are not. The following is proved in [1, Theorem 6].

**THEOREM 25 (An Axiomatization for $\langle\mathbb{R};<,\times\rangle$)** *The theory with the axioms $A_<$, $T_<$, $L_<$ (Theorem 5), $A_\times$, $C_\times$, $U_\times$, $I_\times$, $Z_\times$ (Theorem 22), $O_\times$ (Theorem 16), $D_\times^o$, $R_\times^e$ (Theorem 23) along with the following completely axiomatizes the structures $\langle\mathbb{R};0,1,-1,<,\times,{}^{-1}\rangle$.*

$(P_<)$ $\forall x\,(0 < x \to \exists y \neq 0[x = y^2])$
$(N_<)$ $\exists u\,(-1 < 0 < 1 < u)$
$(O_\times^-)$ $\forall x, y, z\,(z < 0 \wedge x < y \to y \cdot z < x \cdot z)$

The axiomatizability of the structure $\langle\mathbb{Q};<,\times\rangle$ seemed to be missing (or ignored) in the literature. Since $\langle\mathbb{Q};<,+,\times\rangle$ is not decidable (Corollary 21), one could not immediately infer the decidability of $\langle\mathbb{Q};<,\times\rangle$. Also, $+$ is not definable in $\langle\mathbb{Q};<,\times\rangle$, this follows from Theorem 26 below, and so Corollary 21 cannot imply its undecidability. The decidability of $\langle\mathbb{Q};<,\times\rangle$ was proved, and an explicit axiomatization was provided for it, for the first time in [1, Theorem 7]:

**THEOREM 26 (An Axiomatization for $\langle \mathbb{Q}; <, \times \rangle$)** *The theory with $\mathbf{A}_<$, $\mathbf{T}_<$, $\mathbf{L}_<$ (Theorem 5), $\mathbf{O}_\times$ (Theorem 16), $\mathbf{A}_\times$, $\mathbf{C}_\times$, $\mathbf{U}_\times$, $\mathbf{I}_\times$, $\mathbf{Z}_\times$ (Theorem 22), $\mathbf{R}_\times^e$ (Theorem 23), $\mathbf{M}_\times$ (Theorem 24), $\mathbf{N}_<$ (Theorem 25), along with the following completely axiomatizes the structure $\langle \mathbb{Q}; \mathbf{0}, \mathbf{1}, -\mathbf{1}, <, \times, {}^{-1} \rangle$.*

$(\mathbf{D}_<^\times)$  $\{\forall x, y \exists v\,(\mathbf{0} < x < y \rightarrow x < v^n < y)\}_{n>0}$

The axioms $\mathbf{M}_\times$ in Theorem 26 state that $\mathbb{Q}^+$ is dense in the set of its positive radicals.

**THEOREM 27 (Non-Axiomatizability of $\langle \mathbb{N}; <, \times \rangle$, $\langle \mathbb{Z}; <, \times \rangle$)** *The full first-order theory of $\langle \mathbb{N}; <, \times \rangle$ and $\langle \mathbb{Z}; <, \times \rangle$ are not axiomatizable by any decidable set of sentences.*

For a proof note that successor and zero are definable in both of these structures by $v = \mathfrak{s}(u) \iff u < v \wedge \neg \exists w [u < w < v]$ and $(u = \mathbf{0}) \iff u \times \mathfrak{s}(u) = u$, respectively. So, $+$ is definable in $\langle \mathbb{N}; <, \times \rangle$ by Tarski-Robinson's identity [19]: $(u+v=w) \iff [u=v=w=\mathbf{0}] \vee [w \neq \mathbf{0} \wedge \mathfrak{s}(wu)\mathfrak{s}(wv) = \mathfrak{s}(w^2\mathfrak{s}(uv))]$. Thus, by Theorem 20, $\langle \mathbb{N}; <, \times \rangle$ is not axiomatizable; neither is $\langle \mathbb{Z}; <, \times \rangle$ since $\mathbb{N}$ is definable in it by the formula $\mathbf{0} \leqslant v$.

The exponential function is not total in $\mathbb{Z}$ or $\mathbb{Q}$, even when the base is positive: $2^{-1} \notin \mathbb{Z}$ and $2^{(1/2)} \notin \mathbb{Q}$. As for $\mathbb{N}$ we take $\mathfrak{exp}(x, y) = x^y$ with the convention that $0^0 = 1$; and of course $0^x = 0$ for every $x > 0$. For $\mathbb{R}$ and $\mathbb{C}$ we consider $x \mapsto e^x$ for the Napier-Euler number $e$ in the place of $\mathfrak{exp}(x)$, since if $x$ is negative, then the value of $x^y$ may not exist in $\mathbb{R}$, such as $(-1)^{(1/4)}$, and even if it exists in $\mathbb{C}$ it may not be unique (for example, $1^{(1/4)}$ could be $1$, $-1$, $i'$, or $-i'$); indeed one can take any positive real number for $e$. We also add $+$ and $\times$ to the language; so, by *the real exponential field* we mean $\langle \mathbb{R}; +, \times, e^x \rangle$ and by *the complex exponential field* we mean $\langle \mathbb{C}; +, \times, e^x \rangle$.

**THEOREM 28 (Non-Axiomatizability of $\langle \mathbb{N}; \mathfrak{exp} \rangle$)** *The first-order theory of $\langle \mathbb{N}; \mathfrak{exp} \rangle$ is not axiomatizable.*

Since one can define $\times$ and $+$ in the structure $\langle \mathbb{N}; \mathfrak{exp} \rangle$ (see [6, Exercise 1, page 223]) by $(u \times v = w) \iff \forall x [x^w = (x^u)^v]$ and $(u + v = w) \iff \forall x [x^w = (x^u) \times (x^v)]$. So, the result follows from Theorem 20.

**THEOREM 29 (Non-Axiomatizability of $\langle \mathbb{C}; +, \times, e^x \rangle$)** *The complex exponential field is not axiomatizable.*

Indeed, the formula $\forall x, y\,(e^{xy} = -x^2 = 1 \rightarrow e^{xy \cdot v} = 1)$ defines $\mathbb{Z}$ in $\langle \mathbb{C}; +, \times, e^x \rangle$, see e.g. [15], since $\forall x\,(x^2 = -1 \leftrightarrow x = \pm i')$ holds in $\mathbb{C}$, and for every $z$ we have $[e^{\pm i'z} = 1$ if and only if $z = k\pi$ for some $k \in \mathbb{Z}]$. The result follows from Corollary 21.

One of the most exciting questions in the axiomatizability theory is the question of the axiomatizability of $\langle \mathbb{R}; +, \times, e^x \rangle$, the real exponential field (due to Tarski) which is still open. An interesting instance of interaction between seemingly different areas of mathematics (number theory and logic) is the result of Macintyre and Wilkie [14] which states that $\langle \mathbb{R}; +, \times, e^x \rangle$ *is axiomatizable if and only if Weak Schanuel's Conjecture is true*. So, if a computer scientist or a mathematical logician shows the (non-)axiomatizability of $\langle \mathbb{R}; +, \times, e^x \rangle$, then weak Schanuel's conjecture is solved in computational number theory, and if a number theorist solves that problem, then we know whether $\langle \mathbb{R}; +, \times, e^x \rangle$ is axiomatizable or not. If the conjecture is true, then we have an axiomatization for $\langle \mathbb{R}; +, \times, e^x \rangle$ which is "quite complicated and ugly" according to Marker [15].

## § IDENTITIES (OVER $+, \times, \mathfrak{exp}$ IN $\mathbb{R}^+$).

First-order sentences can be restricted in at least two ways: one can consider the sentences of the form

(a) $\exists \vec{x}\, \eta(\vec{x})$ where $\eta(\vec{x})$ is an equation (between two terms on $\vec{x}$ and possibly some other parameters); or

(b) $\forall \vec{x}\, \eta(\vec{x})$ where $\eta(\vec{x})$ is as above.

The (a) formulas are called *diophantine equations* and are closely related to Hilbert's 10th Problem. Since these formulas are discussed elsewhere (see e.g. [2, 5, 18]), here we discuss the formulas in (b), which are called *identities*.

For a proof of the parts (i) and (ii) of the following theorem see e.g. [11]; and for a proof of part (iii), which is due to Martin [16], see e.g. [12, Corollary 3.7].

**THEOREM 30 (Identities With A Single Operation.)**
*(i) The identities of $\langle \mathbb{R}^+; + \rangle$ are axiomatized by*

$(\mathbf{A}_+)$  $x + (y+z) = (x+y) + z$

$(\mathbf{C}_+)$  $x + y = y + x$

*(ii) The identities of $\langle \mathbb{R}^+; \mathbf{1}, \times \rangle$ are axiomatized by*

$(\mathbf{A}_\times)$  $x \cdot (y \cdot z) = (x \cdot y) \cdot z$

$(\mathbf{C}_\times)$  $x \cdot y = y \cdot x$

$(\mathbf{U}_\times)$  $x \cdot \mathbf{1} = x$

*(iii) The identities of $\langle \mathbb{R}^+; \mathbf{1}, \mathfrak{exp} \rangle$ are axiomatized by*

$(\mathbf{C}_\wedge)$  $(x^y)^z = (x^z)^y$

$(\mathbf{Z}_\wedge)$  $\mathbf{1}^x = \mathbf{1}$

$(\mathbf{U}_\wedge)$  $x^\mathbf{1} = x$

Let us note that $0 \notin \mathbb{R}^+$ and so the identity $(\mathbf{U}_+)$ $x + \mathbf{0} = x$ is not expressible here; and since we do not have $-$ in our language, the identity $(\mathbf{I}_+)$ $x + (-x) = \mathbf{0}$ is not expressible either. The part (I) of the following theorem appears in [11]; for the part (II), which appeared in [16] first, see e.g. [12, Corollary 3.9].

**THEOREM 31 (Identities With Two Operations.)**
*(I) The identities of $\langle \mathbb{R}^+; \mathbf{1}, +, \times \rangle$ are axiomatized by $\mathbf{A}_+$, $\mathbf{C}_+$, $\mathbf{A}_\times$, $\mathbf{C}_\times$, $\mathbf{U}_\times$ (Theorem 30) along with the following identity:*

$(\mathbf{D}_+^\times)$  $x \cdot (y+z) = (x \cdot y) + (x \cdot z)$

*(II) The identities of $\langle \mathbb{R}^+; \mathbf{1}, \times, \mathfrak{exp} \rangle$ are axiomatized by $\mathbf{A}_\times$, $\mathbf{C}_\times$, $\mathbf{U}_\times$, $\mathbf{Z}_\wedge$, $\mathbf{U}_\wedge$ (Theorem 30) along with the following identities:*

$(\mathbf{D}_\wedge^\times)$  $x^{(y \cdot z)} = (x^y)^z$

$(\mathbf{D}_\times^\wedge)$  $(x \cdot y)^z = x^z \cdot y^z$

Let us note that the axiom $\mathbf{C}_\wedge$ (Theorem 30.iii) is provable from $\mathbf{C}_\times$ (Theorem 30.ii) and $\mathbf{D}_\wedge^\times$ (Theorem 31.II). The axioms in Theorem 31.I (for $\{+, \times\}$) suffices for proving many of the high-school identities, such as:

- The Binomial Identity: $(x+y)^n = \sum_{i \leqslant n} \binom{n}{i} x^i y^{n-i}$,
- $(x+y+1)^n = \sum_{(i+j \leqslant n)} \binom{n}{i+j} \binom{i+j}{i} x^i y^j$,

and the more difficult one:

$(\mathbf{W}_v^u)$: $(A^u+B^u)^v(C^v+D^v)^u = (A^v+B^v)^u(C^u+D^u)^v$, where $A(x)=x+1$, $B(x)=x^2+x+1$, $C(x)=x^3+1$, and $D(x)=x^4+x^2+1$ are polynomials on $x$.

Of course the $\{+,\times\}$-identities in Theorem 31.I can prove $(\mathbf{W}_v^u)$ when both $u$ and $v$ are positive natural numbers. We now show that the identities of Table 2 derive $(\mathbf{W}_v^u)$ when at least one of $u$ or $v$ is a natural number. So, we assume that say $u \in \mathbb{N}$ and note that $AD = BC = x^5+x^4+x^3+x^2+x+1$. We have:

$(A^u+B^u)^v(C^v+D^v)^u = (A^u+B^u)^v \sum_{i \leqslant u} \binom{u}{i} C^{vi} D^{v(u-i)} =$

$\sum_{i \leqslant u} \binom{u}{i} [(A^u+B^u) C^i D^{u-i}]^v =$

$\sum_{i \leqslant u} \binom{u}{i} ([(AC)^i (AD)^{u-i}] + [(BC)^i (BD)^{u-i}])^v =$

$\sum_{i \leqslant u} \binom{u}{i} ([(AC)^i (BC)^{u-i}] + [(AD)^i (BD)^{u-i}])^v =$

$\sum_{i \leqslant u} \binom{u}{i} (C^u [A^i B^{u-i}] + D^u [A^i B^{u-i}])^v =$

$\sum_{i \leqslant u} \binom{u}{i} ([C^u+D^u][A^i B^{u-i}])^v =$

$(C^u+D^u)^v \sum_{i \leqslant u} \binom{u}{i} (A^v)^i (B^v)^{u-i} = (A^v+B^v)^u (C^u+D^u)^v$.

Indeed, Wilkie's identity $(\mathbf{W}_v^u)$ is true even when both $u,v$ are variables: since for $E(x) = x^2-x+1$ we have $C = AE$ and $D = BE$, thus $E^{uv}$ can be factored out from both sides of $(\mathbf{W}_v^u)$. Note that the positive-valued polynomial $E$ is not expressible in the language $\{\mathbf{1},+,\times,\mathfrak{exp}\}$.

Tarski's High-School Problem asked whether the identities of Table 2 could axiomatize all the identities of the positive cone of the real exponential field $\langle \mathbb{R}^+; \mathbf{1},+,\times,\mathfrak{exp} \rangle$. It was posed first by Doner & Tarski (1969) and was popularized in 1977 by Henkin [11] as a then open problem. Wilkie [24] showed in 1981 that $(\mathbf{W}_v^u)$ is not derivable from Tarski's high-school identities when both $u$ and $v$ are variables (see also [8]).

Wilkie [24] also proved that the identities of $\langle \mathbb{R}^+; \mathbf{1},+,\times,\mathfrak{exp} \rangle$ are axiomatizable by a decidable set of identities and Gurevič [9] showed that it is not axiomatizable by any finite set of identities. However, Tarski's conjecture holds true for a wide range of identities.

Let us say that a term $t$ over $\{\mathbf{1},+,\times,\mathfrak{exp}\}$ is of level 1 when for every sub-term $u^v$ of $t$ either $u$ is a variable or $u$ contains no variable; for example, $x^\alpha + (1+1)^\beta$. A term $t$ is of level 2 when for every sub-term $u^v$ of $t$ we have that $u$ is of level 1; for example the term $p(x)^u + q(x)^v$ is of level 2 when $p,q$ are polynomials of the variable $x$ and $u,v$ are variables. Let us note that the term $(p(x)^u + q(x)^u)^v$, which appears in $(\mathbf{W}_v^u)$, is not of level 2 in general. The following theorem is proved in [10, Proposition 4.4.5]:

**THEOREM 32 (Tarski's Conjecture For Terms of Level 2)**
*If $(r = s)$ is a valid identity of the structure $\langle \mathbb{R}^+; \mathbf{1},+,\times,\mathfrak{exp} \rangle$ where $r$ and $s$ are terms of level 2, then $(r = s)$ can be proved from the identities of Table 2.*

So, Wilkie's result [24] (Theorem 33 below) is a boundary result, since some terms in Wilkie's identity $(\mathbf{W}_v^u)$ are of level 3 (which are the terms with the property that for every sub-term $u^v$ of them, $u$ is a term of level 2).

**THEOREM 33 (Tarski's Conjecture Not for Higher Levels)**
*The identity $(\mathbf{W}_v^u)$ holds in $\langle \mathbb{R}^+; \mathbf{1},+,\times,\mathfrak{exp} \rangle$ but is not provable from the identities of Table 2 when $u,v,x$ are all variables.*

§ **APPENDIX (An Axiomatization for $\langle \mathbb{N}^+; \times \rangle$).**

An axiomatization for $\langle \mathbb{N}^+; \mathbf{1}, \times \rangle$ was presented in [4] whose proofs are available only in French; an English exposition of the axioms without any proof appears in [23, § III.5]. We need the following notation for presenting the axioms:

$y \sqsubseteq x \iff \exists w(y \cdot w = x)$,
$\mathscr{P}(x) \iff x \neq \mathbf{1} \land \forall y (y \sqsubseteq x \to y = \mathbf{1} \lor y = x)$,
$\mathscr{R}(x,y) \iff \mathscr{P}(x) \land x \sqsubseteq y \land \forall z (\mathscr{P}(z) \land z \neq x \to z \not\sqsubseteq y)$, and
$\mathscr{V}(x,y,z) \iff \mathscr{R}(x,z) \land z \sqsubseteq y \land \forall w (\mathscr{R}(x,w) \land w \sqsubseteq y \to w \sqsubseteq z)$;

which state, respectively, that "$y$ divides $x$", "$x$ is a prime", "$y$ is a power of the prime $x$", and "$z$ is the largest power of the prime $x$ that divides $y$". Here are Cégielski's axioms ([4]):

$(\mathbf{A}_\times)$ $\forall x,y,z\,(x \cdot (y \cdot z) = (x \cdot y) \cdot z)$
$(\mathbf{C}_\times)$ $\forall x,y\,(x \cdot y = y \cdot x)$
$(\mathbf{U}_\times)$ $\forall x\,(x \cdot \mathbf{1} = x)$
$(\mathbf{C}^\times)$ $\forall x,y,z\,(x \cdot y = x \cdot z \to y = z)$
$(\mathbf{U}^\times)$ $\forall x,y\,(x \cdot y = \mathbf{1} \to x = y = \mathbf{1})$
$(\mathbf{D}^\times)$ $\{\forall x,y\,(x^n = y^n \to x = y)\}_{n>1}$
$(\mathbf{E}_\times)$ $\{\forall x \exists u,v\,(x = u^n v \land \forall y,z [x = y^n z \to v \sqsubseteq z])\}_{n>1}$
$(\mathbf{P}^\times)$ $\forall x \exists v\,(\mathscr{P}(v) \land v \not\sqsubseteq x)$
$(\mathbf{R}_\times)$ $\forall u,x,y\,(\mathscr{R}(u,x) \land \mathscr{R}(u,y) \to x \sqsubseteq y \lor y \sqsubseteq x)$
$(\mathbf{V}_\exists)$ $\forall u,x\,[\mathscr{P}(u) \to \exists v \mathscr{V}(u,x,v)]$
$(\mathbf{V}_\sqsubseteq)$ $\forall x,y\,(\forall u,v,w[\mathscr{P}(u) \land \mathscr{V}(u,x,v) \land \mathscr{V}(u,y,w) \to v \sqsubseteq w] \longrightarrow x \sqsubseteq y)$
$(\mathbf{V}_\times)$ $\forall x,y\,(\forall u,v,w[\mathscr{P}(u) \land \mathscr{V}(u,x,v) \land \mathscr{V}(u,y,w) \longrightarrow \mathscr{V}(u,x \cdot y, v \cdot w)])$
$(\mathbf{T}_\times)$ $\forall x,y \exists z \forall u\,(\mathscr{P}(u) \longrightarrow [u \not\sqsubseteq x \to \mathscr{V}(u,z,\mathbf{1})] \land [u \sqsubseteq x \to \forall v\{\mathscr{V}(u,z,v) \leftrightarrow \mathscr{V}(u,y,v)\}])$
$(\mathbf{S}_\times)$ $\{\forall x,y \exists z \forall u\,(\mathscr{P}(u) \longrightarrow [u \sqsubseteq x \cdot y \land \exists v,w\{\mathscr{V}(u,x,v) \land \mathscr{V}(u,y,w^n v)\} \to \mathscr{V}(u,z,u)] \land [\neg(u \sqsubseteq x \cdot y \land \exists v,w\{\mathscr{V}(u,x,v) \land \mathscr{V}(u,y,w^n v)\}) \to \mathscr{V}(u,z,\mathbf{1})])\}_{n>0}$

By $\mathbf{U}_\times$, $\mathbf{C}^\times$, and $\mathbf{U}^\times$ the relation $\sqsubseteq$ is anti-symmetric: if $a \sqsubseteq b \sqsubseteq a$, then $a = b$. For every prime $u$ and every $x$ there exists some $v$, by $\mathbf{V}_\exists$, such that $\mathscr{V}(u,x,v)$. That $v$ is unique by $\mathbf{V}_\sqsubseteq$; so let us denote it by $\mathcal{V}(u,x)$. So, if $u$ ranges over the primes, then $x = \prod_{u \sqsubseteq x} \mathcal{V}(u,x)$. Thus, $\mathbf{V}_\times$ is equivalent to $\mathcal{V}(u,xy) = \mathcal{V}(u,x)\mathcal{V}(u,y)$; and the number $z$ in $\mathbf{T}_\times$ is $\prod_{u \sqsubseteq x} \mathcal{V}(u,y)$. The axiom $\mathbf{S}_\times$ states the existence of $\prod_{[u \sqsubseteq xy, \mathcal{V}(u,x) \sqsubseteq_n \mathcal{V}(u,y)]} u$, where $a \sqsubseteq_n b$ is by definition $\exists w(aw^n = b)$. Finally, we note that the following sentences are provable from the axioms:

$(\mathbf{V}_=)$ $\forall x,y\,(\forall u[\mathscr{P}(u) \to \mathcal{V}(u,x) = \mathcal{V}(u,y)] \longrightarrow x = y)$
$(\mathbf{I}_\times)$ $\forall x \exists w \forall u\,(\mathscr{P}(u) \to [u \not\sqsubseteq x \to \mathcal{V}(u,w) = \mathbf{1}] \land [u \sqsubseteq x \to \mathcal{V}(u,w) = u\mathcal{V}(u,x)])$
$(\mathbf{P}^\times_\exists)$ $\forall x\,(x \neq \mathbf{1} \to \exists u[\mathscr{P}(u) \land u \sqsubseteq x])$

In fact, $\mathbf{V}_=$ follows from $\mathbf{V}_\sqsubseteq$, and $\mathbf{I}_\times$ follows from $\mathbf{S}_\times$ by putting $w = xz$ where $z$ is stated to exist by $\mathbf{S}_\times$ for $x = y, n = 1$. Indeed, $\mathbf{V}_=$ is the axiom $A11$ in [4] ($V_2$ in [23]), and $\mathbf{I}_\times$ is the axiom $A15$ in [4] ($I$ in [23]) which, as we saw, are redundant. For $\mathbf{P}^\times_\exists$ we note that if no prime divides $\alpha \neq 1$, then $\mathcal{V}(u,\alpha) = 1$ for every prime $u$; so by $\mathbf{V}_\sqsubseteq$ we have $\alpha \sqsubseteq y$ for every $y$, and this contradicts $\mathbf{P}^\times$ (by which there are infinitely many primes).

---

**Acknowledgements:** This research was partially supported by a grant from $\mathbb{IPM}$ (No. 98030022).

SAEED SALEHI

Research Institute for Fundamental Sciences (RIFS),
University of Tabriz, 29 Bahman Boulevard,
P.O.Box 51666–17766, Tabriz, Iran.
School of Mathematics,
Institute for Research in Fundamental Sciences ($\mathbb{IPM}$),
P.O.Box 19395–5746, Tehran, Iran.
root@saeedsalehi.ir
http://saeedsalehi.ir/